# RENORMALIZED SELF-INTERSECTION LOCAL TIME FOR FRACTIONAL BROWNIAN MOTION

By Yaozhong Hu[1] and David Nualart[2]

*University of Kansas and Universitat de Barcelona*

Let $B_t^H$ be a $d$-dimensional fractional Brownian motion with Hurst parameter $H \in (0,1)$. Assume $d \geq 2$. We prove that the renormalized self-intersection local time

$$\ell = \int_0^T \int_0^t \delta(B_t^H - B_s^H)\,ds\,dt - \mathbb{E}\Big(\int_0^T \int_0^t \delta(B_t^H - B_s^H)\,ds\,dt\Big)$$

exists in $L^2$ if and only if $H < 3/(2d)$, which generalizes the Varadhan renormalization theorem to any dimension and with any Hurst parameter. Motivated by a result of Yor, we show that in the case $3/4 > H \geq \frac{3}{2d}$, $r(\varepsilon)\ell_\varepsilon$ converges in distribution to a normal law $N(0,T\sigma^2)$, as $\varepsilon$ tends to zero, where $\ell_\varepsilon$ is an approximation of $\ell$, defined through (2), and $r(\varepsilon) = |\log \varepsilon|^{-1}$ if $H = 3/(2d)$, and $r(\varepsilon) = \varepsilon^{d-3/(2H)}$ if $3/(2d) < H$.

**1. Introduction.** The fractional Brownian motion on $\mathbb{R}^d$ with Hurst parameter $H \in (0,1)$ is a $d$-dimensional Gaussian process $B^H = \{B_t^H, t \geq 0\}$ with zero mean and covariance function given by

$$\mathbb{E}(B_t^{H,i} B_s^{H,j}) = \frac{\delta_{ij}}{2}(t^{2H} + s^{2H} - |t-s|^{2H}),$$

where $i,j = 1,\ldots,d$, and $s,t \geq 0$. We will assume that $d \geq 2$. The *self-intersection local time* of $B^H$ is formally defined as

(1) $$I = \int_0^T \int_0^t \delta_0(B_t^H - B_s^H)\,ds\,dt,$$

Received November 2003; revised May 2004.
[1]Supported in part by NSF Grant DMS-02-04613.
[2]Supported in part by MCyT Grant BFM2000-0598. This work was carried out at the University of Kansas.
*AMS 2000 subject classifications.* 60G15, 60G18, 60F05, 60F25, 60H30.
*Key words and phrases.* Fractional Brownian motion, self-intersection local time, Wiener chaos development, renormalization, central limit theorem.

where $\delta_0(x)$ is the Dirac delta function. It measures the amount of time that the process spends intersecting itself on the time interval $[0,T]$ and has been an important topic of the theory of stochastic process. A rigorous definition of this random variable may be obtained by approximating the Dirac function by the heat kernel

$$p_\varepsilon(x) = (2\pi\varepsilon)^{-d/2}\exp\{-|x|^2/(2\varepsilon)\},$$

as $\varepsilon > 0$ tends to zero. We denote the approximated self-intersection local time by

$$(2) \qquad I_\varepsilon = \int_0^T \int_0^t p_\varepsilon(B_t^H - B_s^H)\,ds\,dt,$$

and a natural question is to study the behavior of $I_\varepsilon$ as $\varepsilon$ tends to zero.

For $H = 1/2$, the process $B^H$ is a classical Brownian motion. The self-intersection local time of the Brownian motion has been studied by many authors [see Albeverio, Hu and Zhou (1997), Calais and Yor (1987), He, Yang, Yao and Wang (1995), Hu (1996), Imkeller, Pérez-Abreu and Vives (1995), Varadhan (1969), Yor (1985) and the references therein]. In the case of the planar Brownian motion, Varadhan (1969) has proved that $I_\varepsilon$ does not converge in $L^2$ but it can be renormalized so that $I_\varepsilon - \frac{T}{2\pi}\log(1/\varepsilon)$ converges in $L^2$ as $\varepsilon$ tends to zero. The limit is called the *renormalized self-intersection local time* of the planar Brownian motion. This result has been extended by Rosen (1987) to the (planar) fractional Brownian motion, where it is proved that for $1/2 < H < 3/4$, $I_\varepsilon - C_H T \varepsilon^{-1+1/(2H)}$ converges in $L^2$ as $\varepsilon$ tends to zero, where $C_H$ is a constant depending only on $H$ [see also (4): $C_H = C_{H,2}$].

On the other hand, for classical Brownian motion and when the dimension $d \geq 3$, Yor (1985) and Calais and Yor (1987) have proved that the random variables

$$\begin{cases} (\log(1/\varepsilon))^{-1/2}(I_\varepsilon - \mathbb{E}(I_\varepsilon)), & \text{if } d = 3, \\ \varepsilon^{d/2-3/2}(I_\varepsilon - \mathbb{E}(I_\varepsilon)), & \text{if } d > 3, \end{cases}$$

converge in law to normal distribution (with a different approximation).

The first aim of this paper is to extend Rosen's and Varadhan's results to arbitrary dimensions and with Hurst parameter $H < 3/4$. Along this line, we have established the following result.

THEOREM 1. *Let $I_\varepsilon$ be the random variable defined in* (2). *We have:*

(i) *If $H < 1/d$, then $I_\varepsilon$ converges in $L^2$ as $\varepsilon$ tends to zero.*
(ii) *If $1/d < H < 3/(2d)$, then*

$$(3) \qquad I_\varepsilon - TC_{H,d}\varepsilon^{-d/2+1/(2H)}$$

*converges in $L^2$ as $\varepsilon$ tends to zero, where*

$$(4) \qquad C_{H,d} = (2\pi)^{-d/2}\int_0^\infty (z^{2H}+1)^{-d/2}\,dz.$$

(iii) *If $1/d = H < 3/(2d)$, then*

$$(5) \qquad I_\varepsilon - \frac{T}{2H(2\pi)^{d/2}} \log(1/\varepsilon)$$

*converges in $L^2$ as $\varepsilon$ tends to zero.*

(iv) *If $H \geq 3/(2d)$, then the difference $I_\varepsilon - \mathbb{E}(I_\varepsilon)$ does not converge in $L^2$.*

That means, if $H < 3/(2d)$, the difference $I_\varepsilon - \mathbb{E}(I_\varepsilon)$ converges in $L^2$ as $\varepsilon$ tends to zero to the *renormalized self-intersection local time*.

In the case $H \geq 3/(2d)$ we have established the following version of the central limit theorem.

THEOREM 2. *Suppose $3/(2d) \leq H < 3/4$. Then the random variables*

$$(6) \qquad \begin{cases} (\log(1/\varepsilon))^{-1/2}(I_\varepsilon - \mathbb{E}(I_\varepsilon)), & \text{if } H = 3/(2d), \\ \varepsilon^{d/2 - 3/(4H)}(I_\varepsilon - \mathbb{E}(I_\varepsilon)), & \text{if } H > 3/(2d), \end{cases}$$

*converge as $\varepsilon$ tends to zero in distribution to a normal law $N(0, T\sigma^2)$, where $\sigma^2$ is a constant depending on $d$ and $H$.*

We conjecture that, as functions of the terminal time $T$, these processes converge in law to $\{\sigma B_T, T \geq 0\}$, where $B_T$ is a Brownian motion. We plan to discuss this problem in a forthcoming paper.

For $H \geq 3/4$ our method of proof of Theorem 2 does not work and the study of the fluctuations of $I_\varepsilon - \mathbb{E}(I_\varepsilon)$ in this case is an open problem.

The proof of Theorem 1 is based on some estimates deduced by Hu in (2001) from the local nondeterminism property of the fractional Brownian motion. The proof of Theorem 2 is more involved, and the main ingredient is to show that the projections on each Wiener chaos converge in law to independent Gaussian random variables. The proof of these convergences is based on a recent general criterion for the convergence in distribution to a normal law of a sequence of multiple stochastic integrals established by Nualart and Peccati (2005) and by Peccati and Tudor (2003), in the multidimensional case. We have extended their results to a slightly more general setting applicable to the renormalization of self-intersection local time of fractional Brownian motion and to this end we have established a general central limit theorem which has its own interest.

In the case of the classical Brownian motion, the convergence in law of the Wiener chaos projections of the random variables appearing in (6) has been proved in de Faria, Drumond and Streit (2000) using a martingale approach which cannot be extended to the case $H \neq 1/2$. The corresponding convergence in law to a family of independent Brownian motions has been established in Rezgui and Streit (2002).

The paper is organized as follows. In Section 2 we prove Theorem 1. Section 3 is devoted to a general central limit theorem using chaos expansions. Some examples are given and will be used to deal with the self-intersection local time of the fractional Brownian motion. Section 4 deduces the chaos expansion of the approximated self-intersection local time, and in Section 5 we prove Theorem 2. The proofs of the main results are based on a sequence of technical lemmas, stated and proved in Section 6.

**2. Renormalized self-intersection local time, case $dH < 3/2$.** Let $B^H = \{B_t^H, t \geq 0\}$ be a $d$-dimensional fractional Brownian motion of Hurst parameter $H \in (0,1)$.

Consider the approximation $I_\varepsilon$ of the self-intersection local time introduced in (2). From the equality

$$p_\varepsilon(x) = (2\pi)^{-d} \int_{\mathbb{R}^d} \exp\{i\langle \xi, x\rangle\} \exp\left\{-\frac{\varepsilon|\xi|^2}{2}\right\} d\xi$$

and the definition of $I_\varepsilon$, we obtain

(7) $\quad I_\varepsilon = (2\pi)^{-d} \int_0^T \int_0^t \int_{\mathbb{R}^d} \exp\{i\langle \xi, B_t^H - B_s^H\rangle\} \exp\left\{-\frac{\varepsilon|\xi|^2}{2}\right\} d\xi \, ds \, dt.$

Therefore,

(8) $\quad \mathbb{E}(I_\varepsilon^2) = (2\pi)^{-2d} \int_{\mathcal{T}} \int_{\mathbb{R}^{2d}} \mathbb{E}(\exp\{i\xi(B_t^H - B_s^H)\} \exp\{i\eta(B_{t'}^H - B_{s'}^H)\})$
$\times \exp\{-(\varepsilon|\xi|^2 + \varepsilon|\eta|^2)/2\} \, d\xi \, d\eta \, ds \, dt \, ds' \, dt',$

where

(9) $\quad \mathcal{T} = \{(s,t,s',t'): 0 < s < t < T, 0 < s' < t' < T\}.$

Throughout this paper we will make use of the following notation, for any $\tau = (s,t,s',t')$:

(10) $\quad \lambda(\tau) = |t-s|^{2H}, \qquad \rho(\tau) = |t'-s'|^{2H}$

and

(11) $\quad \mu(\tau) = \frac{1}{2}[|s-t'|^{2H} + |s'-t|^{2H} - |t-t'|^{2H} - |s-s'|^{2H}].$

Notice that $\lambda$ is the variance of $B_t^{H,1} - B_s^{H,1}$, $\rho$ is the variance of $B_{t'}^{H,1} - B_{s'}^{H,1}$ and $\mu$ is the covariance between $B_t^{H,1} - B_s^{H,1}$ and $B_{t'}^{H,1} - B_{s'}^{H,1}$, where $B^{H,1}$ denotes a one-dimensional fractional Brownian motion with Hurst parameter $H$.

With this notation, for any $\xi, \eta \in \mathbb{R}^d$, we can write

(12) $\quad \mathbb{E}[(\langle \xi, B_t^H - B_s^H\rangle + \langle \eta, B_{t'}^H - B_{s'}^H\rangle)^2] = \lambda|\xi|^2 + \rho|\eta|^2 + 2\mu\langle \xi, \eta\rangle.$

As a consequence, from (8) and (12) we deduce for all $\varepsilon > 0$,

$$
\begin{aligned}
\mathbb{E}(I_\varepsilon^2) &= (2\pi)^{-2d} \int_\mathcal{T} \int_{\mathbb{R}^{2d}} e^{-((\lambda+\varepsilon)|\xi|^2 + 2\mu\langle\xi,\eta\rangle + (\rho+\varepsilon)|\eta|^2)/2} \, d\xi \, d\eta \, d\tau \\
&= (2\pi)^{-d} \int_\mathcal{T} [(\lambda+\varepsilon)(\rho+\varepsilon) - \mu^2]^{-d/2} \, d\tau.
\end{aligned}
\tag{13}
$$

On the other hand, the expectation of the random variable $I_\varepsilon$ is given by

$$
\begin{aligned}
\mathbb{E}(I_\varepsilon) &= \int_0^T \int_0^t p_{\varepsilon+|t-s|^{2H}}(0) \, ds \, dt \\
&= (2\pi)^{-d/2} \int_0^T \int_0^t (\varepsilon + |t-s|^{2H})^{-d/2} \, ds \, dt \\
&= (2\pi)^{-d/2} \int_0^T (T-s)(\varepsilon + s^{2H})^{-d/2} \, ds.
\end{aligned}
\tag{14}
$$

PROOF OF THEOREM 1. Assertion (i) follows easily from (13) and (14). From (14), making the change of variables $s = z\varepsilon^{1/(2H)}$, we obtain, if $1/d < H < 3/(2d)$,

$$
\begin{aligned}
\mathbb{E}(I_\varepsilon) &= \frac{\varepsilon^{1/(2H)-d/2}}{(2\pi)^{d/2}} \int_0^{T\varepsilon^{1/(2H)}} (T - z\varepsilon^{1/(2H)})(z^{2H}+1)^{-d/2} \, dz \\
&= \varepsilon^{1/(2H)-d/2} T C_{H,d} + o(\varepsilon).
\end{aligned}
$$

For $H = 1/d$ we get

$$
\mathbb{E}(I_\varepsilon) = \frac{T\log(1/\varepsilon)}{2H(2\pi)^{d/2}} + o(\varepsilon).
$$

Hence, the convergence in $L^2$ of the random variables (3) and (5) is equivalent to the convergence of $I_\varepsilon - \mathbb{E}(I_\varepsilon)$.

From (13) and (14) we obtain

$$
\mathbb{E}(I_\varepsilon I_\eta) - \mathbb{E}(I_\varepsilon)\mathbb{E}(I_\eta) \\
= (2\pi)^{-d} \int_\mathcal{T} [((\lambda+\varepsilon)(\rho+\eta) - \mu^2)^{-d/2} - ((\lambda+\varepsilon)(\rho+\eta))^{-d/2}] \, d\tau.
$$

Therefore, a necessary and sufficient condition for the convergence in $L^2$ of $I_\varepsilon - \mathbb{E}(I_\varepsilon)$ is that

$$
\Xi_T =: \int_\mathcal{T} [(\lambda\rho - \mu^2)^{-d/2} - (\lambda\rho)^{-d/2}] \, d\tau < \infty.
\tag{15}
$$

Finally, Lemma 11 allows us to conclude the proof of the theorem. $\square$

**3. Central limit theorem via chaos expansion.** In this section we will establish a general central limit theorem for nonlinear functionals of a Gaussian process, based on the Wiener chaos expansion. The proof of this theorem uses some recent results by Nualart and Peccati (2005) and by Peccati and Tudor (2003).

Consider a separable Hilbert space $H$. Let $\{e_k, k \geq 1\}$ be a complete orthonormal system in $H$. For every $n \geq 1$, we denote by $H^{\odot n}$ the $n$th symmetric tensor product of $H$. For every $p = 0, \ldots, n$ and for every $f \in H^{\odot n}$, we define the *contraction* of $f$ of order $p$ to be the element of $H^{\otimes 2(n-p)}$ defined by

$$f \otimes_p f = \sum_{i_1, \ldots, i_p = 1}^{\infty} \langle f, e_{i_1} \otimes \cdots \otimes e_{i_p} \rangle_{H^{\otimes p}} \otimes \langle f, e_{i_1} \otimes \cdots \otimes e_{i_p} \rangle_{H^{\otimes p}}.$$

Suppose that $X = \{X(h), h \in H\}$ is an *isonormal Gaussian process* on $H$. This means that $X$ is a centered Gaussian family indexed by the elements of $H$, defined on some probability space $(\Omega, \mathcal{F}, \mathbb{P})$ and such that, for every $h, h' \in H$,

$$\mathbb{E}(X(h)X(h')) = \langle h, h' \rangle_H.$$

Assume that $\mathcal{F}$ is a $\sigma$-algebra generated by $X$.

For every $n \geq 1$ we will denote by $I_n$ the multiple Wiener integral [Hu and Kallianpur (1998)].

THEOREM 3. *Consider a sequence of square integrable and centered random variables $\{F_k, k \geq 1\}$ with Wiener chaos expansions*

$$F_k = \sum_{n=1}^{\infty} I_n(f_{n,k}).$$

*Suppose that:*

(i) $\lim_{N \to \infty} \limsup_{k \to \infty} \sum_{n=N+1}^{\infty} n! \|f_{n,k}\|^2_{H^{\otimes n}} = 0$;
(ii) *for every $n \geq 1$*, $\lim_{k \to \infty} n! \|f_{n,k}\|^2_{H^{\otimes n}} = \sigma_n^2$;
(iii) $\sum_{n=1}^{\infty} \sigma_n^2 = \sigma^2 < \infty$;
(iv) *for all $n \geq 2$, $p = 1, \ldots, n-1$*, $\lim_{k \to \infty} \|f_{n,k} \otimes_p f_{n,k}\|^2_{H^{\otimes 2(n-p)}} = 0$.

*Then $F_k$ converges in distribution to the normal law $N(0, \sigma^2)$ as $k$ tends to infinity.*

PROOF. By Theorem 1 in Nualart and Peccati (2005), conditions (ii) and (iv) imply that for each fixed $n \geq 1$ the sequence of multiple integrals $I_n(f_{n,k})$ converges in distribution to the normal law $N(0, \sigma_n^2)$ as $k$ tends to infinity. Furthermore, by Theorem 1 in Peccati and Tudor (2003), for each $n \geq 1$, the vector $(I_1(f_{1,k}), \ldots, I_n(f_{n,k}))$ converges in law to the vector $(\xi_1, \ldots, \xi_n)$, where $\{\xi_n, n \geq 1\}$ are independent centered Gaussian random variables with variances $\{\sigma_n^2, n \geq 1\}$. Now let $f$ be a $C^1$ function such that $|f|$ and $|f'|$ are bounded by 1. For every $N \geq 1$, set

$$F_k^{(N)} = \sum_{n=1}^{N} I_n(f_{n,k}).$$

Set $\xi^{(N)} = \sum_{n=1}^{N} \xi_n$, and $\xi = \sum_{n=1}^{\infty} \xi_n$. Then

$$|\mathbb{E}(f(F_k)) - \mathbb{E}(f(\xi))|$$
$$\leq |\mathbb{E}(f(F_k)) - \mathbb{E}(f(F_k^{(N)}))|$$
$$+ |\mathbb{E}(f(F_k^{(N)})) - \mathbb{E}(f(\xi^{(N)}))| + |\mathbb{E}(f(\xi^{(N)})) - \mathbb{E}(f(\xi))|$$
$$\leq \left(\sum_{n=N+1}^{\infty} n!\|f_{n,k}\|_{H^{\otimes n}}^2\right)^{1/2} + |\mathbb{E}(f(F_k^{(N)})) - \mathbb{E}(f(\xi^{(N)}))|$$
$$+ |\mathbb{E}(f(\xi^{(N)})) - \mathbb{E}(f(\xi))|.$$

Then, taking first the limit as $k$ tends to infinity, and then the limit as $N$ tends to infinity, and applying conditions (i), (iii) and the convergence in law of the vector $(I_1(f_{1,k}), \ldots, I_n(f_{n,k}))$, we get the desired result. □

REMARK 1. A sufficient condition for (i) and (iii) is that for all $n \geq 1$, and $k \geq 1$,

$$n!\|f_{n,k}\|_{H^{\otimes n}}^2 \leq \delta_n,$$

where $\sum_{n=1}^{\infty} \delta_n < \infty$.

REMARK 2. The assumptions of Theorem 3 imply the convergence in law of the projections on the Wiener chaos to a family of independent Gaussian random variables.

The following propositions contain examples of applications of the above central limit theorem which will be useful in proving the renormalization result for the fractional Brownian motion when $dH \geq 3/2$.

Let $\mathcal{H}$ be the Hilbert space defined as the closure of set $\mathcal{E}$ of step functions on $\mathbb{R}_+$ with respect to the scalar product

$$\langle \mathbf{1}_{[0,t]}, \mathbf{1}_{[0,s]}\rangle_{\mathcal{H}} = \tfrac{1}{2}(t^{2H} + s^{2H} - |t-s|^{2H}).$$

Let $B^H = \{B_t^H, t \geq 0\}$ be a $(2m)$-dimensional fractional Brownian motion with Hurst parameter $H$, with $m \geq 1$. Then, $B^H$ is an isonormal Gaussian process associated with the Hilbert space $H = \mathcal{H}^{2m}$. We denote by $I_n$ the multiple stochastic integral which provides an isometry between the symmetric tensor product $(\mathcal{H}^{2m})^{\odot}$ equipped with the norm $\sqrt{n!}\|\cdot\|_{H^{\otimes n}}$ and the $n$th Wiener chaos of $B^H$.

We will make use of the following notation:
$$K_1(x, y, z) = |z + y|^{2H} + |z - x|^{2H} - |z + y - x|^{2H} - z^{2H},$$
$$K_2(x, z) = K_1(x, x, z).$$

PROPOSITION 4. *Let $\{B_t^H, t \geq 0\}$ be a $(2m)$-dimensional fractional Brownian motion with Hurst parameter $H \in (0, 3/4)$. For any $0 < s < t$, we set*

$$\text{(16)} \qquad Z_{s,t} = \prod_{i=1}^{2m}(B_t^{H,i} - B_s^{H,i}).$$

*Define, for any fixed $x > 0$,*
$$\Lambda_T(x) = \frac{1}{\sqrt{T}}\int_0^T Z_{s,s+x}\,ds.$$

*Then as $T$ tends to infinity $\Lambda_T(x)$ converges in distribution to a normal law $N(0, G(x, x))$, where*
$$G(x, x) = 2^{1-2m}\int_0^\infty K_2^{2m}(x, y)\,dy.$$

PROOF. We have
$$\mathbb{E}[\Lambda_T^2(x)] = \frac{2^{1-2m}}{T}\int_0^T\int_0^t K_2^{2m}(x, t-s)\,ds\,dt$$
$$= 2^{1-2m}\int_0^T K_2^{2m}(x, y)(1 - y/T)\,dy$$

and, clearly, this converges to $G(x, x)$ as $T$ tends to infinity. We can write
$$\Lambda_T(x) = I_{2m}(f_{2m,T}),$$
where
$$f_{2m,T}(\mathbf{i}_{2m}, r_1, \ldots, r_{2m}) = \frac{1}{(2m)!\sqrt{T}}\int_0^T \prod_{i=1}^{2m}\mathbf{1}_{[s,s+x]}(r_i)\,ds,$$

if all components of the multi-index $\mathbf{i}_{2m} = (i_1, \ldots, i_{2m})$ are different and zero otherwise. Then, by Theorem 1 of Nualart and Peccati (2005), it suffices to show that for any $p$, $1 \leq p \leq 2m - 1$,

$$\text{(17)} \qquad \lim_{T \to \infty}\|f_{2m,T}\otimes_p f_{2m,T}\|_{H^{\otimes(2m-p)}}^2 = 0.$$

The contraction of order $p$ of the kernel $f_{2m,T}$ is given by

$$(f_{2m,T} \otimes_p f_{2m,T})(\mathbf{i}_{2m-p}, r_1, \ldots, r_{2m-p}, \mathbf{j}_{2m-p}, u_1, \ldots, u_{2m-p})$$

$$= \frac{p!}{((2m)!)^2 T 2^p} \int_0^T \int_0^T K_2^p(x, t-s) \prod_{j=1}^{2m-p} \mathbf{1}_{[s,s+x]}(r_j) \mathbf{1}_{[t,t+x]}(u_j) \, ds \, dt.$$

In order to show (17) it suffices to prove that

(18)
$$\lim_{T \to \infty} \frac{1}{T} \int_{0 < t_1 < t_2 < t_3 < t_4 < T} K_2^p(x, t_2 - t_1) K_2^p(x, t_4 - t_3)$$
$$\times K_2^{2m-p}(x, t_3 - t_2) K_2^{2m-p}(x, t_4 - t_1) \, dt = 0.$$

By L'Hôpital's rule, (18) holds if

(19)
$$\lim_{T \to \infty} \frac{1}{T} \int_{0 < t_1 < t_2 < t_3 < T} K_2^p(x, t_2 - t_1) K_2^p(x, T - t_3)$$
$$\times K_2^{2m-p}(x, t_3 - t_2) K_2^{2m-p}(x, T - t_1) \, dt_1 \, dt_2 \, dt_3$$
$$= 0.$$

Set $t_2 - t_1 = a$, $t_3 - t_2 = b$, $T - t_3 = c$. In order to show (19) it suffices to show that

(20)
$$\lim_{T \to \infty} \frac{1}{T} \int_{0 < a+b+c < T} K_2^p(x, a) K_2^p(x, c) K_2^{2m-p}(x, b)$$
$$\times K_2^{2m-p}(x, a+b+c) \, da \, db \, dc = 0.$$

To show (20) we will make use of the following inequalities:

(21) $$|K_2(x, a)| \leq k_H a^{2H-2} x^2,$$

if $x < \frac{a}{2}$, where $k_H = 2^{3-2H} H(2H-1)$,

(22) $$a + b + c \geq 3(abc)^{1/3}$$

and

(23) $$|K_2(x, a)| \leq 2(ax)^H.$$

We decompose the set $G = \{0 < a + b + c < T\}$ as follows:

$$G = \bigcup_{i=1}^{8} G_i,$$

where

$$G_1 = \{(a,b,c) \in G : a > 2x, b > 2x, c > 2x\},$$
$$G_2 = \{(a,b,c) \in G : a > 2x, b > 2x, c < 2x\},$$
$$G_3 = \{(a,b,c) \in G : a > 2x, b < 2x, c > 2x\},$$
$$G_4 = \{(a,b,c) \in G : a < 2x, b > 2x, c > 2x\},$$
$$G_5 = \{(a,b,c) \in G : a > 2x, b < 2x, c < 2x\},$$
$$G_6 = \{(a,b,c) \in G : a < 2x, b > 2x, c < 2x\},$$
$$G_7 = \{(a,b,c) \in G : a < 2x, b < 2x, c > 2x\},$$
$$G_8 = \{(a,b,c) \in G : a < 2x, b < 2x, c < 2x\}.$$

Set, for $i = 1, \ldots, 8$,

$$\Lambda_T^i = \frac{1}{T} \int_{G_i} K_2^p(x,a) K_2^p(x,c) K_2^{2m-p}(x,b) K_2^{2m-p}(x,a+b+c) \, da \, db \, dc.$$

In what follows, we will denote by $k$ a generic constant that may depend on $H$ and $x$.

*Case* 1. Consider first the integral over the region $G_1$. Then, using (21) and (22) we obtain, for $T$ large enough,

$$\Lambda_T^1 \leq \frac{k}{T} \int_{G_1} a^{(2H-2)p} b^{(2H-2)p} c^{(2H-2)(2m-p)} (a+b+c)^{(2H-2)(2m-p)} \, da \, db \, dc$$
$$\leq \frac{k}{T} \left( \int_{2x}^T a^{4/3(2H-2)p} \, da \right)^2 \left( \int_{2x}^T a^{4/3(2H-2)(2m-p)} \, da \right)$$
$$\leq k T^{4/3(2H-2)(2m+p)+2},$$

which converges to zero as $T$ tends to infinity, because

$$\tfrac{4}{3}(2H-2)(2m+p) + 2 \leq \tfrac{4}{3}(2H-2)3 + 2 = 8H - 6 < 0.$$

*Case* 2. Using (21), (22) and (23) we obtain

$$\Lambda_T^2 \leq \frac{k}{T} \int_{G_2} a^{(2H-2)p} b^{(2H-2)p} c^{H(2m-p)} (a+b+c)^{(2H-2)(2m-p)} \, da \, db \, dc$$
$$\leq k T^{(4H-4)(m+p/2)+1},$$

which converges to zero as $T$ tends to infinity, because

$$(4H-4)(m+p/2) + 1 \leq 3(2H-2) + 1 = 6H - 5 < 0.$$

*Case* 3. Using (21), (22) and (23) we obtain

$$\Lambda_T^3 \le \frac{k}{T}\int_{G_3} a^{(2H-2)p}b^{Hp}c^{(2H-2)(2m-p)}(a+b+c)^{(2H-2)(2m-p)}\,da\,db\,dc$$

$$\le kT^{(4H-4)(2m-p/2)+1} \to 0,$$

because $2m - p/2 \ge 2m - (2m-1)/2 = m + 1/2 \ge 1$, and, so

$$(4H-4)(2m-p/2) + 1 \le 4H - 3 < 0.$$

*Case* 4. The integral over the set $G_4$ is handled as in Case 3.

*Case* 5. Using (21), (22) and (23) we obtain

$$\Lambda_T^5 \le \frac{k}{T}\int_{G_5} a^{(2H-2)p}b^{Hp}c^{H(2m-p)}(a+b+c)^{(2H-2)(2m-p)}\,da\,db\,dc$$

$$= kT^{(4H-4)m} \to 0.$$

Finally, the integrals over the sets $G_6$, $G_7$ and $G_8$ are treated in a similar way.

□

PROPOSITION 5. *Consider the processes $\Lambda_T(x)$ defined in Proposition* 4. *If $3/(2d) < H < 3/4$, then the random variables*

$$Y_T := \int_0^\infty (1+x^{2H})^{-d/2-m}\Lambda_T(x)\,dx$$

*converge in law to the normal distribution $N(0,\sigma^2)$, as $T$ tends to infinity, where*

$$(24) \quad \sigma^2 := \int_0^\infty \int_0^\infty (1+x^{2H})^{-d/2-m}(1+y^{2H})^{-d/2-m}G(x,y)\,dx\,dy$$

*and*

$$G(x,y) := 2^{1-2m}\int_0^\infty K_1^{2m}(x,y,z)\,dz.$$

PROOF. Notice that

$$\mathbb{E}[\Lambda_T(x)\Lambda_T(y)] = \frac{2^{1-2m}}{T}\int_0^T\int_0^t K_1^{2m}(x,y,t-s)\,ds\,dt$$

$$= 2^{1-2m}\int_0^T K_1^{2m}(x,y,z)(1-z/T)\,dz,$$

and this clearly converges as $T$ tends to infinity to $G(x,y)$.

As a consequence, Proposition 4 and Proposition 2 of Peccati and Tudor (2003) imply that the finite-dimensional distributions of the process $\{\Lambda_T(x), x > 0\}$ converge to the finite-dimensional distributions of a Gaussian centered process with covariance function $G(x, y)$. In order to complete the proof of the lemma we need to show the following facts:

(i)
$$\int_0^\infty (1 + x^{2H})^{-d/2-m} \sup_T \mathbb{E}(|\Lambda_T(x)|) < \infty. \tag{25}$$

(ii) For all $K > 0$,
$$\lim_{T \to \infty} \mathbb{E}[(\Lambda_T(x) - \Lambda_T(y))^2] \to 0 \tag{26}$$

as $|x - y| < \delta \to 0$, and $x, y \leq K$.

(iii) $\sigma^2 < \infty$.

PROOF OF (i). Property (25) follows from
$$\int_0^\infty (1 + x^{2H})^{-d/2-m} G(x,x)^{1/2}\, dx$$
$$\leq \int_0^\infty (1 + x^{2H})^{-d/2-m} \left(\int_0^{2x} K_2^{2m}(x,y)\, dy\right)^{1/2} dx$$
$$+ \int_0^\infty (1 + x^{2H})^{-d/2-m} \left(\int_{2x}^\infty K_2^{2m}(x,y)\, dy\right)^{1/2} dx \tag{27}$$
$$\leq k_{H,m} \int_0^\infty (1 + x^{2H})^{-d/2-m} \left(\int_0^{2x} (y^{4Hm} + x^{4Hm})\, dy\right)^{1/2} dx$$
$$+ k_{H,m} \int_0^\infty (1 + x^{2H})^{-d/2-m} \left(\int_{2x}^\infty y^{(2H-2)2m} x^{4m}\, dy\right)^{1/2} dx$$
$$\leq k'_{H,m} \int_0^\infty (1 + x^{2H})^{-d/2-m} x^{2mH+1/2}\, dx < \infty,$$

because $Hd > 3/2$. □

PROOF OF (ii). We have
$$\sup_T \mathbb{E}[(\Lambda_T(x) - \Lambda_T(y))^2]$$
$$\leq 2^{1-2m} \int_0^\infty |K_2^{2m}(x,z) + K_2^{2m}(y,z) - 2K_1^{2m}(x,y,z)|\, dz,$$

and, by dominated convergence this tends to zero as $|x - y|$ tends to zero in the compact interval $[0, K]$. □

PROOF OF (iii). This follows from (27) and the fact that, by Cauchy–Schwarz inequality, we have $G(x,y) \le \sqrt{G(x,x)G(y,y)}$.  □

Let us now show that (i)–(iii) imply the convergence in law of $Y_T$ to $\xi$, where $\xi$ is a random variable with the $N(0,\sigma^2)$ distribution. Let $f$ be a function such that $|f|$ and $|f'|$ are bounded by 1. We make the decomposition

$$|\mathbb{E}(f(Y_T)) - \mathbb{E}(f(\xi))|$$
$$\le |\mathbb{E}(f(Y_T)) - \mathbb{E}(f(Y_T^{(K)}))| + |\mathbb{E}(f(Y_T^{(K)})) - \mathbb{E}(f(\xi^{(K)}))|$$
$$+ |\mathbb{E}(f(\xi^{(K)})) - \mathbb{E}(f(\xi))|,$$

where

$$Y_T^{(K)} = \int_0^K (1+x^{2H})^{-d/2-m} \Lambda_T(x)\, dx,$$

$\xi^{(K)}$ has the law $N(0,\sigma_K^2)$, and

$$\sigma_K^2 = 2^{1-2m} \int_0^K \int_0^K (1+x^{2H})^{-d/2-m}(1+y^{2H})^{-d/2-m} G(x,y)\, dx\, dy.$$

Property (i) implies that $\lim_{K\to\infty} \sup_T |\mathbb{E}(f(Y_T)) - \mathbb{E}(f(Y_T^{(K)}))| = 0$.

Property (iii) implies that $\lim_{K\to\infty} |\mathbb{E}(f(\xi^{(K)})) - \mathbb{E}(f(\xi))| = 0$, and property (ii) implies that for any fixed $K>0$, $\lim_{T\to\infty} |\mathbb{E}(f(Y_T^{(K)})) - \mathbb{E}(f(\xi^{(K)}))| = 0$.

This completes the proof of the proposition.  □

PROPOSITION 6.  *Assume $H = 3/(2d) < 3/4$. Consider the processes*

$$X_T = \frac{1}{\sqrt{T\log T}} \int_0^\infty (1+x^{2H})^{-d/2-m} \left( \int_0^T Z_{s,s+x}\, ds \right) dx,$$

*where $T>0$, and $Z_{s,s+x}$ has been defined in* (16). *Then $X_T$ converges in distribution, as $T$ tends to infinity, to a normal law $N(0,\sigma^2)$, where*

$$\sigma^2 := 2^{1-2m} \int_0^\infty \int_0^\infty (xy)^{-3/2-2mH} K_1^{2m}(x,y,1)\, dx\, dy.$$

PROOF.  As in the proof of Proposition 4 we can show that, for each fixed $x>0$, the random variables $\frac{1}{\sqrt{T\log T}} \int_0^T Z_{s,s+x}\, ds$ converge to 0 in $L^2$ as $T$ tends to infinity. For this reason we cannot apply the same method of proof as in Proposition 4, and new ideas are required. The basic ingredient of the proof will be a scaling argument and Theorem 1 of Nualart and Peccati (2005).

Let us see first that the asymptotic behavior of $X_T$ as $T$ tends to infinity is the same as $Y_T$, where

$$Y_T := \frac{1}{\sqrt{T \log T}} \int_1^\infty x^{-3/2 - 2mH} \left( \int_0^T Z_{s,s+x} \, ds \right) dx.$$

In fact,

$$\mathbb{E}[|X_T - Y_T|] \leq \frac{k}{\sqrt{T \log T}} \int_0^\infty (1 + x^{2H})^{-5/2 - m} \mathbb{E}\left( \left| \int_0^T Z_{s,s+x} \, ds \right| \right) dx,$$

and using the estimate (27) it is not difficult to see that this converges to zero as $T$ tends to infinity.

So, it suffices to show that $Y_T$ converges in distribution as $T$ tends to infinity to the normal law $N(0, \sigma^2)$. The proof will be done in several steps.

*Step* 1. We claim that

(28) $$\lim_{T \to \infty} \mathbb{E}(Y_T^2) = \sigma^2 < \infty,$$

where $\sigma^2$ is defined in (24). In fact,

$$\mathbb{E}(Y_T^2) = \frac{2}{T \log T} \int_1^\infty \int_1^\infty (xy)^{-3/2 - 2mH} \int_0^T \int_0^t \mu(s, s+x, t, t+y) \, ds \, dt \, dx \, dy$$

$$= \frac{2^{1-2m}}{\log T} \int_1^\infty \int_1^\infty (xy)^{-3/2 - 2mH} \int_0^T K_1^{2m}(x, y, z) \left( 1 - \frac{z}{T} \right) dz \, dx \, dy.$$

Thus, by L'Hôpital's rule

$$\lim_{T \to \infty} \mathbb{E}(Y_T^2) = \lim_{T \to \infty} T 2^{1-2m} \int_1^\infty \int_1^\infty (xy)^{-3/2 - 2mH} K_1^{2m}(x, y, T) \, dx \, dy$$

$$= \lim_{T \to \infty} 2^{1-2m} \int_{1/T}^\infty \int_{1/T}^\infty (xy)^{-3/2 - 2mH} K_1^{2m}(x, y, 1) \, dx \, dy$$

$$= \sigma^2.$$

*Step* 2. Consider the decomposition

(29) $$Y_T = Y_{T,1} + Y_{T,2},$$

where

$$Y_{T,1} = \frac{1}{\sqrt{T \log T}} \int_1^\infty x^{-3/2 - 2mH} \left( \int_0^{T/2} Z_{s,s+x} \, ds \right) dx,$$

$$Y_{T,2} = \frac{1}{\sqrt{T \log T}} \int_1^\infty x^{-3/2 - 2mH} \left( \int_{T/2}^T Z_{s,s+x} \, ds \right) dx.$$

From the scaling and the stationary increments properties of fractional Brownian motion, it follows that the random variables $Y_{T,1}$ and $Y_{T,2}$ have the same distribution as $Y_{T,3}$, where

$$Y_{T,3} = \frac{1}{\sqrt{2T \log T}} \int_2^\infty x^{-3/2 - 2mH} \left( \int_0^T Z_{s,s+x} \, ds \right) dx.$$

We are going to compute the second and fourth moments of $Y_T$, using the decomposition (29). We have

$$\mathbb{E}(Y_T^2) = \mathbb{E}(Y_{T,1}^2) + \mathbb{E}(Y_{T,2}^2) + 2\mathbb{E}(Y_{T,1}Y_{T,2})$$
$$= 2\mathbb{E}(Y_{T,3}^2) + 2\mathbb{E}(Y_{T,1}Y_{T,2}).$$

Clearly, $\lim_{T \to \infty} \mathbb{E}(Y_{T,3}^2) = \sigma^2/2$, and this implies that

(30) $$\lim_{T \to \infty} \mathbb{E}(Y_{T,1}Y_{T,2}) = 0.$$

Consider now the fourth moment of $Y_T$:

(31) $$\mathbb{E}(Y_T^4) = 2\mathbb{E}(Y_{T,3}^4) + 4\mathbb{E}(Y_{T,1}^3 Y_{T,2}) + 4\mathbb{E}(Y_{T,1} Y_{T,2}^3) + 6\mathbb{E}(Y_{T,1}^2 Y_{T,2}^2).$$

*Step* 3. We claim that

(32) $$\lim_{T \to \infty} \mathbb{E}(Y_{T,1}^2 Y_{T,2}^2) = \sigma^4/4,$$
$$\lim_{T \to \infty} \mathbb{E}(Y_{T,1}^3 Y_{T,2}) = 0, \qquad \lim_{T \to \infty} \mathbb{E}(Y_{T,1} Y_{T,2}^3) = 0.$$

Let us prove the first identity of (32). The random variables $(Y_{T,1}Y_{T,2})$ have the same joint distribution as $(Z_{T,1}Z_{T,2})$, where

$$Z_{T,1} = \frac{1}{\sqrt{\log T}} \int_{1/T}^\infty x^{-3/2 - 2mH} \left( \int_0^{1/2} Z_{s,s+x} \, ds \right) dx,$$

$$Z_{T,2} = \frac{1}{\sqrt{\log T}} \int_{1/T}^\infty x^{-3/2 - 2mH} \left( \int_{1/2}^1 Z_{s,s+x} \, ds \right) dx.$$

As a consequence, we have

(33) $$\mathbb{E}(Y_{T,1}^2 Y_{T,2}^2) = (\log T)^{-2} \int_\mathcal{D} \mathbb{E}\left( \prod_{i=1}^4 \Delta_i \right) d\mu,$$

where $\Delta_i = B_{t_i + x_i} - B_{t_i}$, for $i = 1, 2, 3, 4$, $\mathcal{D} = [1/T, \infty)^4 \times [0, 1/2]^2 \times [1/2, 1]^2$ and $d\mu = \prod_{i=1}^4 x_i^{-3/2 - 2mH} \, dx \, dt$. Using the notation $d_{ij} = \mathbb{E}(\Delta_i \Delta_j)$ for $i, j = 1, \ldots, 4$, we can write

$$\mathbb{E}\left( \prod_{i=1}^4 \Delta_i \right) = d_{12}d_{34} + d_{13}d_{24} + d_{14}d_{23}$$

and

$$(34) \qquad \left(\mathbb{E}\left(\prod_{i=1}^{4}\Delta_i\right)\right)^{2m} = (d_{12}d_{34})^{2m} + R.$$

We have

$$(\log T)^{-2}\int_{\mathcal{D}}(d_{12}d_{34})^{2m}\,d\mu = \mathbb{E}(Y_{T,1}^2)\mathbb{E}(Y_{T,2}^2),$$

and this converges to $\sigma^4/4$ as $T$ tends to infinity. On the other hand, the residual term $R$ in (34) does not contribute to the limit (33). In fact, by Hölder's inequality we obtain, for any indexes $(i,j,k)$ such that $i+j+k = 2m$, and $j \geq 1$ or $k \geq 1$,

$$(\log T)^{-2}\left|\int_{\mathcal{D}}(d_{12}d_{34})^i(d_{13}d_{24})^j(d_{14}d_{23})^k\,d\mu\right|$$
$$\leq (\log T)^{-2}\left(\int_{\mathcal{D}}(d_{12}d_{34})^{2m}\,d\mu\right)^{i/(2m)}$$
$$\times \left(\int_{\mathcal{D}}(d_{13}d_{24})^{2m}\,d\mu\right)^{j/(2m)}\left(\int_{\mathcal{D}}(d_{14}d_{23})^{2m}\,d\mu\right)^{k/(2m)}.$$

Then, this converges to zero as $T$ tends to infinity because of (30), $j+k \geq 1$ and

$$(\log T)^{-2}\int_{\mathcal{D}}(d_{14}d_{23})^{2m}\,d\mu = (\log T)^{-2}\int_{\mathcal{D}}(d_{13}d_{24})^{2m}\,d\mu = (\mathbb{E}(Y_{T,1}Y_{T,2}))^2.$$

This completes the proof of the first identity of (32). The proof of the other two identities of (32) follows the same argument.

*Step* 4. By the same arguments as above and Minkowski inequality it follows that

$$\mathbb{E}(Y_T^4) = (\log T)^{-2}\int_{\mathcal{E}}(d_{12}d_{34} + d_{13}d_{24} + d_{14}d_{23})^{2m}\,d\mu$$
$$\leq 3^{2m}(\log T)^{-2}\int_{\mathcal{E}}(d_{12}d_{34})^{2m}\,d\mu$$
$$= 3^{2m}(\mathbb{E}(Y_T^2))^2,$$

where $\mathcal{E} = [1/T,\infty)^4 \times [0,1]^4$. So, the sequences $Y_T$ and $Y_{T,3}$ are bounded in $L^4$, and

$$(35) \qquad \mathbb{E}((Y_T - Y_{T,3})^4) \leq 3^{2m}(\mathbb{E}(Y_T - Y_{T,3})^2)^2,$$

which clearly converges to zero as $T$ tends to infinity. Then, from (35), (31) and (32) we obtain $\lim_{T\to\infty}\mathbb{E}(Y_T^4) = 3\sigma^4$. Finally, the result follows from Theorem 1 of Nualart and Peccati (2005). $\square$

**4. Wiener chaos expansion of the self-intersection local time.** In this section we will first compute the Wiener chaos expansion of the approximation of the self-intersection local time $I_\varepsilon$ defined by (2).

Given a multi-index $\mathbf{i}_n = (i_1, \ldots, i_n)$, $1 \leq i_j \leq d$, we set

$$\alpha(\mathbf{i}_n) = \mathbb{E}[X_{i_1} \cdots X_{i_n}],$$

where the $X_i$ are independent $N(0,1)$ random variables. Notice that

$$\alpha(\mathbf{i}_{2m}) = \frac{(2m_1)! \cdots (2m_d)!}{(m_1)! \cdots (m_d)! 2^m},$$

if $n = 2m$ is even and for each $k = 1, \ldots, d$, the number of components of $\mathbf{i}_{2m}$ equal to $k$, denoted by $2m_k$, is also even, and $\alpha(\mathbf{i}_n) = 0$, otherwise.

LEMMA 7. *We have*

$$I_\varepsilon = \sum_{m=0}^{\infty} I_{2m}(f_{2m,\varepsilon}),$$

*where $f_{2m,\varepsilon}$ is the element of $(\mathcal{H}^d)^{\otimes 2m}$ given by*

(36)
$$\begin{aligned}
&f_{2m,\varepsilon}(\mathbf{i}_{2m}, r_1, \ldots, r_{2m}) \\
&= \frac{(2\pi)^{-d/2} \alpha(\mathbf{i}_{2m})}{(2m)!} \int_0^T \int_0^t ds\, dt (\varepsilon + |t-s|^{2H})^{-d/2-m} \prod_{j=1}^{2m} \mathbf{1}_{[s,t]}(r_j).
\end{aligned}$$

PROOF. Stroock's formula leads to

$$f_{n,\varepsilon} = \frac{1}{n!} \int_0^T \int_0^t \mathbb{E}[D_{r_1,\ldots,r_n}^{i_1,\ldots,i_n} p_\varepsilon(B_t^H - B_s^H)]\, ds\, dt,$$

where $i_j \in \{1, \ldots, d\}$, and $r_l \in [0, T]$.

Let us compute this expectation

$$\mathbb{E}[D_{r_1,\ldots,r_n}^{i_1,\ldots,i_n} p_\varepsilon(B_t^H - B_s^H)] = \mathbb{E}[\partial^{i_1} \cdots \partial^{i_n} p_\varepsilon(B_t^H - B_s^H)] \prod_{j=1}^n \mathbf{1}_{[s,t]}(r_j).$$

On the other hand, using the Fourier transform

$$\begin{aligned}
&\mathbb{E}[\partial^{i_1} \cdots \partial^{i_n} p_\varepsilon(B_t^H - B_s^H)] \\
&= \frac{i^n}{(2\pi)^d} \int_{\mathbb{R}^d} \xi^{i_1} \cdots \xi^{i_n} \mathbb{E}[e^{i\langle \xi, B_t^H - B_s^H \rangle}] e^{-\varepsilon |\xi|^2/2}\, d\xi \\
&= \frac{i^n}{(2\pi)^d} \int_{\mathbb{R}^d} \xi^{i_1} \cdots \xi^{i_n} \exp\left\{-\frac{1}{2}(|t-s|^{2H} + \varepsilon)|\xi|^2\right\} d\xi \\
&= i^n (2\pi)^{-d/2} (\varepsilon + |t-s|^{2H})^{-d/2-n/2} \mathbb{E}[X_{i_1} \cdots X_{i_n}],
\end{aligned}$$

and the result follows. □

Let us compute the $L^2$-norm of the $n$th Wiener chaos of $I_\varepsilon$.

$$\mathbb{E}[(I_{2m}(f_{2m,\varepsilon}))^2] = (2m)!\|f_{2m,\varepsilon}\|^2_{H^{\otimes(2m)}}$$

(37)
$$= (2m)! \sum_{m_1+\cdots+m_d=m} \frac{(2m)!}{(2m_1)!\cdots(2m_d)!} \frac{(2\pi)^{-d}}{((2m)!)^2} \alpha(\mathbf{i}_{2m})^2$$

$$\times \int_{\mathcal{T}} (\varepsilon+\lambda)^{-d/2-m} (\varepsilon+\rho)^{-d/2-m} \mu^{2m}\, d\tau$$

$$= \frac{\alpha_m}{(2\pi)^d 2^{2m}} \int_{\mathcal{T}} (\varepsilon+\lambda)^{-d/2-m} (\varepsilon+\rho)^{-d/2-m} \mu^{2m}\, d\tau,$$

where

$$\alpha_m = \sum_{m_1+\cdots+m_d=m} \frac{(2m_1)!\cdots(2m_d)!}{(m_1!)^2\cdots(m_d!)^2}.$$

If $Hd < 3/2$, we can deduce the Wiener chaos expansion of the renormalized self-intersection local time $\ell$:

$$\ell = \sum_{m=1}^\infty I_{2m}(f_{2m}),$$

where

$$f_{2m}(\mathbf{i}_{2m}, r_1, \ldots, r_{2m})$$
$$= \frac{(2\pi)^{-d/2}\alpha(\mathbf{i}_{2m})}{(2m)!} \int_0^T \int_0^t ds\, dt\, |t-s|^{-Hd-2mH} \prod_{j=1}^{2m} \mathbf{1}_{[s,t]}(r_j).$$

**5. Proof of Theorem 2.** The purpose of this section is to show Theorem 2.

CASE $3/(2d) \leq H < 3/4$. We are going to check that the family of random variables $\varepsilon^{d/2-3/(4H)}(I_\varepsilon - \mathbb{E}(I_\varepsilon))$ satisfies the conditions of Theorem 3, for any sequence $\varepsilon(k) \to 0$. We have

$$\varepsilon^{d/2-3/(4H)}(I_\varepsilon - \mathbb{E}(I_\varepsilon)) = \varepsilon^{d/2-3/(4H)} \sum_{m=1}^\infty I_{2m}(f_{2m,\varepsilon}).$$

Thus, from (37) we obtain

$$\varepsilon^{d-3/(2H)} \mathbb{E}[(I_{2m}(f_{2m,\varepsilon}))^2]$$
$$= \frac{\alpha_m \varepsilon^{d-3/(2H)}}{(2\pi)^d 2^{2m}} \int_{\mathcal{T}} [(\varepsilon+\lambda)(\varepsilon+\rho)]^{-d/2-m} \mu^{2m}\, d\tau.$$

Making the change of variables $\tau \to \varepsilon^{1/(2H)}\tau$ yields

$$\varepsilon^{d-3/(2H)}\mathbb{E}[(I_{2m}(f_{2m,\varepsilon}))^2]$$
$$= \frac{\alpha_m \varepsilon^{1/(2H)}}{(2\pi)^d 2^{2m}} \int_{\mathcal{T}^\varepsilon} [(1+\lambda)(1+\rho)]^{-d/2-m} \mu^{2m} \, d\tau,$$

where $\mathcal{T}^\varepsilon = \varepsilon^{-1/(2H)}\mathcal{T}$. Using the decomposition $\mathcal{T}^\varepsilon \cap \{s < s'\} = \mathcal{T}_1^\varepsilon \cup \mathcal{T}_2^\varepsilon \cup \mathcal{T}_3^\varepsilon$, where $\mathcal{T}_i^\varepsilon = \varepsilon^{-1/(2H)}\mathcal{T}_i$ and the sets $\mathcal{T}_i$ are introduced in (46), we obtain

$$\varepsilon^{d-3/(2H)}\mathbb{E}[(I_{2m}(f_{2m,\varepsilon}))^2]$$
$$= \frac{\alpha_m}{(2\pi)^d 2^{2m-1}} \int_{0<a+b+c<T\varepsilon^{-1/(2H)}} (T - (a+b+c)\varepsilon^{1/(2H)}) \Psi_m \, da \, db \, dc,$$

where

(38) $$\Psi_m = \sum_{i=1}^{3} [(1+\lambda_i)(1+\rho_i)]^{-d/2-m} \mu_i^{2m},$$

and the functions $\lambda_i$, $\rho_i$, $\mu_i$, have been defined in (47)–(52).

As a consequence,

(39) $$\lim_{\varepsilon\downarrow 0} \varepsilon^{d-3/(2H)}\mathbb{E}[(I_{2m}(f_{2m,\varepsilon}))^2] = \frac{\alpha_m T}{(2\pi)^d 2^{2m-1}} \int_{\mathbb{R}_+^3} \Psi_m \, da \, db \, dc,$$

and this implies condition (ii) of Theorem 3. We also have

$$\varepsilon^{d-3/(2H)}\mathbb{E}[(I_{2m}(f_{2m,\varepsilon}))^2] \leq \frac{\alpha_m T}{(2\pi)^d 2^{2m-1}} \int_{\mathbb{R}_+^3} \Psi_m \, da \, db \, dc,$$

and, taking into account Remark 1 after Theorem 3, in order to check conditions (i) and (iii) it suffices to show that

$$\sigma^2 := \sum_{m=1}^{\infty} \frac{\alpha_m}{(2\pi)^d 2^{2m-1}} \int_{\mathbb{R}_+^3} \Psi_m \, da \, db \, dc < \infty.$$

From (38) we have

$$\sigma^2 = \sum_{i=1}^{3} \sum_{m=1}^{\infty} \frac{\alpha_m}{(2\pi)^d 2^{2m-1}} \int_{\mathbb{R}_+^3} (\widehat{\lambda}_i \widehat{\rho}_i)^{-d/2} \widehat{\gamma}_i^m \, da \, db \, dc,$$

where $\widehat{\lambda}_i = 1+\lambda_i$, $\widehat{\rho}_i = 1+\rho_i$ and

$$\widehat{\gamma}_i = \frac{\mu_i^2}{(1+\lambda_i)(1+\rho_i)}.$$

Hence,

$$\sigma^2 = \sum_{i=1}^{3} \frac{2}{(2\pi)^d} \int_{\mathbb{R}_+^3} (\widehat{\lambda}_i \widehat{\rho}_i)^{-d/2} \left[ \left( \sum_{m=0}^{\infty} \frac{(2m)! \widehat{\gamma}_i^m}{(m!)^2 2^{2m}} \right)^d - 1 \right] da\, db\, dc.$$

Using Newton's binomial formula and the definitions (64) and (65) we obtain

$$\sigma^2 = \frac{2}{(2\pi)^d} \sum_{i=1}^{3} \int_{\mathbb{R}_+^3} (\widehat{\lambda}_i \widehat{\rho}_i)^{-d/2} [(1 - \widehat{\gamma}_i)^{-d/2} - 1]\, da\, db\, dc$$

(40)
$$= \frac{2}{(2\pi)^d} \sum_{i=1}^{3} \int_{\mathbb{R}_+^3} ((\widehat{\delta}_i)^{-d/2} - (\widehat{\lambda}_i \widehat{\rho}_i)^{-d/2})\, da\, db\, dc$$

$$= \frac{2}{(2\pi)^d} \sum_{i=1}^{3} \int_{\mathbb{R}_+^3} \widehat{\Theta}_i\, da\, db\, dc,$$

where the functions $\widehat{\Theta}_i$ are defined in (65). Finally, $\sigma^2 < \infty$ by Lemma 13.

Condition (iv) of Theorem 3 follows from the convergence (43) in Proposition 8.

CASE $H = 3/(2d)$. As in the proof of the previous case we need to check that the family of random variables $(\log(1/\varepsilon))^{-1/2}(I_\varepsilon - \mathbb{E}(I_\varepsilon))$ satisfies the conditions of Theorem 3, for any sequence $\varepsilon(k) \to 0$. We have

$$(\log(1/\varepsilon))^{-1/2}(I_\varepsilon - \mathbb{E}(I_\varepsilon)) = (\log(1/\varepsilon))^{-1/2} \sum_{m=1}^{\infty} I_{2m}(f_{2m,\varepsilon}),$$

and as before we obtain

$$\frac{\mathbb{E}[(I_{2m}(f_{2m,\varepsilon}))^2]}{\log(1/\varepsilon)}$$
$$= \frac{2\alpha_m}{(2\pi)^d 2^{2m} \log(1/\varepsilon)}$$
$$\quad \times \int_{0 < a+b+c < T\varepsilon^{-1/(2H)}} (T - t(a+b+c)\varepsilon^{-1/(2H)}) \Psi_m\, da\, db\, dc,$$

where $\Psi_m$ is defined by (38). From the estimates obtained in the proof of Lemma 13 it follows that the term $(a+b+c)\varepsilon^{-1/(2H)}$ gives no contribution to the limit and we can write, using L'Hôpital's rule,

$$\lim_{\varepsilon \downarrow 0} \frac{\mathbb{E}[(I_{2m}(f_{2m,\varepsilon}))^2]}{\log(1/\varepsilon)}$$
$$= \lim_{\varepsilon \downarrow 0} \frac{2T\alpha_m}{(2\pi)^d 2^{2m} \log(1/\varepsilon)} \int_{0 < a+b+c < T\varepsilon^{-1/(2H)}} \Psi_m\, da\, db\, dc$$

$$= \lim_{x\uparrow\infty} \frac{2T\alpha_m}{(2\pi)^d 2^{2m} \log x} \int_{\substack{0<a+b<\eta \\ \eta<Tx^{1/(2H)}}} \Psi_m(a,b,\eta-a-b)\,da\,db\,d\eta$$

$$= \lim_{x\uparrow\infty} \frac{x^{1/(2H)}T^2\alpha_m}{H2^{2m}(2\pi)^d} \int_{0<a+b<Tx^{1/(2H)}} \Psi_m(a,b,Tx^{1/(2H)}-a-b)\,da\,db$$

$$= \lim_{x\uparrow\infty} \frac{x^{3/(2H)}T^4\alpha_m}{H2^{2m}(2\pi)^d} \int_{0<\alpha+\beta<1} \Psi_m(Tx^{1/(2H)}(\alpha,\beta,1-\alpha-\beta))\,d\alpha\,d\beta.$$

Set $e_{\alpha,\beta} = (\alpha, \beta, 1-\alpha-\beta)$. We have $3/(2H) = d$, and

$$x^d \Psi_m(Tx^{1/(2H)} e_{\alpha,\beta})$$

$$= T^{-3} \sum_{i=1}^{3} ((T^{-2H}x^{-1} + \lambda_i)(T^{-2H}x^{-1} + \rho_i))^{-d/2-m} \mu_i^{2m}(e_{\alpha,\beta}),$$

and when $x$ tends to infinity this converges to

$$T^{-3} \sum_{i=1}^{3} (\lambda_i \rho_i)^{-d/2} \gamma_i^m(e_{\alpha,\beta}).$$

Hence,

(41)
$$\lim_{\varepsilon\downarrow 0} \frac{\mathbb{E}[(I_{2m}(f_{2m,\varepsilon}))^2]}{\log(1/\varepsilon)}$$
$$= \frac{T\alpha_m}{H2^{2m}(2\pi)^d} \sum_{i=1}^{3} \int_{0<\alpha+\beta<1} (\lambda_i \rho_i)^{-d/2} \gamma_i^m(e_{\alpha,\beta})\,d\alpha\,d\beta.$$

This proves condition (ii) of Theorem 3.

On the other hand,

$$\sigma^2 := \sum_{i=1}^{3} \sum_{m=1}^{\infty} \frac{\alpha_m}{H2^{2m}(2\pi)^d} \int_{0<\alpha+\beta<1} (\lambda_i \rho_i)^{-d/2} \gamma_i^m(e_{\alpha,\beta})\,d\alpha\,d\beta$$

$$= \sum_{i=1}^{3} \frac{1}{H(2\pi)^d} \int_{0<\alpha+\beta<1} (\lambda_i \rho_i)^{-d/2} \left[ \left( \sum_{m=0}^{\infty} \frac{(2m)!\gamma_i^m}{(m!)^2 2^{2m}} \right)^d - 1 \right](e_{\alpha,\beta})\,d\alpha\,d\beta.$$

Using Newton's binomial formula and the definitions (64) and (65) we obtain

$$\sigma^2 = \frac{1}{H(2\pi)^d} \sum_{i=1}^{3} \int_{0<\alpha+\beta<1} (\lambda_i \rho_i)^{-d/2} [(1-\gamma_i)^{-d/2} - 1](e_{\alpha,\beta})\,d\alpha\,d\beta$$

(42)
$$= \frac{1}{H(2\pi)^d} \sum_{i=1}^{3} \int_{0<\alpha+\beta<1} ((\delta_i)^{-d/2} - (\lambda_i \rho_i)^{-d/2})(e_{\alpha,\beta})\,d\alpha\,d\beta$$

$$= \frac{1}{H(2\pi)^d} \sum_{i=1}^{3} \int_{0<\alpha+\beta<1} \Theta_i(e_{\alpha,\beta}) \, d\alpha \, d\beta.$$

Finally, $\sigma^2 < \infty$ by Lemma 12. This proves condition (iii).

Finally,

$$\lim_{\varepsilon \downarrow 0} \frac{1}{\log(1/\varepsilon)} \sum_{m=N+1}^{\infty} \mathbb{E}[(I_{2m}(f_{2m,\varepsilon}))^2]$$

$$\leq \sum_{m=N+1}^{\infty} \sum_{i=1}^{3} \frac{T\alpha_m}{H2^{2m}(2\pi)^d} \int_{0<\alpha+\beta<1} (\widehat{\lambda}_i \widehat{\rho}_i)^{-d/2} \widehat{\gamma}_i^m(e_{\alpha,\beta}) \, d\alpha \, d\beta,$$

and this converges to zero as $N$ tends to infinity. So condition (i) holds. Condition (iv) follows from the convergence (44) in Proposition 8.

PROPOSITION 8. *Fix $m \geq 1$ and $1 \leq p \leq 2m-1$.*

(i) *If $d > 3/(2H)$, then*

(43) $$\lim_{\varepsilon \downarrow 0} \varepsilon^{2d-3/H} \|f_{2m,\varepsilon} \otimes_p f_{2m,\varepsilon}\|^2_{H^{\otimes 2(2m-p)}} = 0;$$

(ii) *if $d = 3/(2H)$ and $H < 3/4$, then*

(44) $$\lim_{\varepsilon \downarrow 0} (\log(1/\varepsilon))^{-2} \|f_{2m,\varepsilon} \otimes_p f_{2m,\varepsilon}\|^2_{H^{\otimes 2(2m-p)}} = 0.$$

PROOF OF (43). The proof will be done in several steps.

*Step* 1. Let us first compute $f_{2m,\varepsilon} \otimes_p f_{2m,\varepsilon}$. We have, from (36),

$$(f_{2m,\varepsilon} \otimes_p f_{2m,\varepsilon})(\mathbf{i}_{2m-p}, \mathbf{j}_{2m-p}, r_1, \ldots, r_{2m-p}, u_1, \ldots, u_{2m-p})$$

$$= \frac{(2\pi)^{-d}}{((2m)!)^2} \sum_{\mathbf{k}_p} \alpha(\mathbf{i}_{2m-p}, \mathbf{k}_p) \alpha(\mathbf{k}_p, \mathbf{j}_{2m-p})$$

$$\times \int_{\mathcal{T}} (\varepsilon + |t-s|^{2H})^{-d/2-m} (\varepsilon + |t'-s'|^{2H})^{-d/2-m} \mu^p(s,t,s',t')$$

$$\times \prod_{j=1}^{2m-p} \mathbf{1}_{[s,t]}(r_j) \mathbf{1}_{[s',t']}(u_j) \, d\tau.$$

As a consequence,

$$\|f_{2m,\varepsilon} \otimes_p f_{2m,\varepsilon}\|^2_{H^{\otimes 2(2m-p)}}$$

(45) $$= A^2 \int_{\mathcal{D}} \prod_{i=1}^{4} (\varepsilon + |t_i - s_i|^{2H})^{-d/2-m}$$

$$\times \mu^p_{1,2}(s,t) \mu^p_{3,4}(s,t) \mu^{2m-p}_{1,3}(s,t) \mu^{2m-p}_{2,4}(s,t) \, ds \, dt,$$

where

$$A^2 = \frac{(2\pi)^{-2d}}{((2m)!)^4} \sum_{\mathbf{i}_{2m-p},\mathbf{j}_{2m-p}} \left(\sum_{\mathbf{k}_p} \alpha(\mathbf{i}_{2m-p},\mathbf{k}_p)\alpha(\mathbf{k}_p,\mathbf{j}_{2m-p})\right)^2,$$

$$s = (s_1, s_2, s_3, s_4), \qquad t = (t_1, t_2, t_3, t_4),$$

$$\mathcal{D} = \{(s,t) \in \mathbb{R}^8 : 0 < s_i < t_i < T\}$$

and

$$\mu_{i,j}(s,t) = \mu(s_i, t_i, s_j, t_j).$$

Set

$$\Phi_\varepsilon = \varepsilon^{2d-3/H} \int_{\mathcal{D}} \prod_{i=1}^{4} (\varepsilon + |t_i - s_i|^{2H})^{-d/2-m}$$
$$\times \mu_{1,2}^p(s,t)\mu_{3,4}^p(s,t)\mu_{1,3}^{2m-p}(s,t)\mu_{2,4}^{2m-p}(s,t)\,ds\,dt.$$

To show (43) it suffices to show that $\lim_{\varepsilon \to 0} \Phi_\varepsilon = 0$.

*Step* 2. We are going to use a representation of $\Phi_\varepsilon$ in terms of the expectation of the product of four random variables. Consider a family of independent one-dimensional fractional Brownian motions with Hurst parameter $H$, $B_t^{i,j}$, $i = 1, 2, 3, 4$, $j = 1, \ldots, 2m$. We can write

$$\Phi_\varepsilon = \varepsilon^{2d-3/H} \mathbb{E}\left(\prod_{i,j} X_\varepsilon^{i,j}\right),$$

where the product is over the pairs $(i,j) = (1,3), (1,4), (2,3)$ and $(3,4)$, and

$$X_\varepsilon^{i,j} = \int_0^T \int_0^t (\varepsilon + |t-s|^{2H})^{-d/2-m} \prod_{\ell=1}^{p} \prod_{k=1}^{2m-p} (B_t^{i,\ell} - B_s^{i,\ell})(B_t^{j,k} - B_s^{j,k})\,ds\,dt.$$

*Step* 3. The variances of the random variables $\varepsilon^{d/2-3/(4H)} X_\varepsilon^{i,j}$ converge to a constant times $T$ as $\varepsilon$ tends to zero. In fact,

$$\mathbb{E}[\varepsilon^{d-3/(2H)}(X_\varepsilon^{i,j})^2] = \frac{(2\pi)^d 2^{2m}}{\alpha_m} \varepsilon^{d-3/(2H)} \mathbb{E}[(I_{2m}(f_{2m,\varepsilon}))^2],$$

and we have already proved in (39) that this converges to

$$2T \int_{\mathbb{R}_+^3} \Psi_m \, da\, db\, dc < \infty.$$

*Step* 4. In order to show that $\lim_{\varepsilon \to 0} \Phi_\varepsilon = 0$, it suffices to show that each family of random variables $\varepsilon^{d/2-3/(4H)} X_\varepsilon^i$ converges in law to a normal random variable. In fact, taking into account that $\mathbb{E}(X_\varepsilon^{i,j} X_\varepsilon^{i',j'}) = 0$, for $(i,j) \neq (i',j')$, and the results proved in Peccati and Tudor (2003), the four-dimensional vector

$$\varepsilon^{d/2-3/(4H)}(X_\varepsilon^{1,3}, X_\varepsilon^{1,4}, X_\varepsilon^{2,3}, X_\varepsilon^{3,4})$$

will converge to a Gaussian random vector with independent components and, therefore, the expectation of the product of its components will converge to zero.

Notice that the families of random variables $\{\varepsilon^{d/2-3/(4H)} X_\varepsilon^{i,j}, \varepsilon > 0\}$ are identically distributed and they have the same distribution as $\{Z_\varepsilon, \varepsilon > 0\}$, where

$$Z_\varepsilon = \varepsilon^{d/2-3/(4H)} \int_0^T \int_0^t (\varepsilon + |t-s|^{2H})^{-d/2-m} \prod_{i=1}^{2m}(B_t^{H,i} - B_s^{H,i}) \, ds \, dt,$$

and $\{B_t^H, t \geq 0\}$ is a $(2m)$-dimensional fractional Brownian motion with Hurst parameter $H$. We make the decomposition

$$Z_\varepsilon = Z_\varepsilon^1 - Z_\varepsilon^2,$$

where

$$Z_\varepsilon^2 = \varepsilon^{1/(4H)} \int_{s < T\varepsilon^{-1/(2H)} < t} (1 + |t-s|^{2H})^{-d/2-m} Z_{s,t} \, ds \, dt,$$

$$Z_\varepsilon^1 = \varepsilon^{1/(4H)} \int_{s < t \wedge T\varepsilon^{-1/(2H)}} (1 + |t-s|^{2H})^{-d/2-m} Z_{s,t} \, ds \, dt,$$

and $Z_{s,t}$ has been defined in (16). Let us show that $Z_\varepsilon^2$ tends to zero in $L^2(\Omega)$. We have

$$\|Z_\varepsilon^2\|_2 \leq \varepsilon^{1/(4H)} \int_{s < T\varepsilon^{-1/(2H)} < t} (1 + |t-s|^{2H})^{-d/2-m} \|Z_{s,t}\|_2 \, ds \, dt$$

$$= \varepsilon^{1/(4H)} \int_{s < T\varepsilon^{-1/(2H)} < t} (1 + |t-s|^{2H})^{-d/2-m} |t-s|^{2mH} \, ds \, dt$$

$$\leq \varepsilon^{1/(4H)} \int_{s < T\varepsilon^{-1/(2H)} < t} (1 + |t-s|^{2H})^{-d/2} \, ds \, dt$$

$$= \varepsilon^{1/(4H)} \int_0^{T\varepsilon^{-1/(2H)}} \int_{T\varepsilon^{-1/(2H)}-s}^\infty (1 + x^{2H})^{-d/2} \, dx \, ds$$

$$= \varepsilon^{1/(4H)} \int_0^{T\varepsilon^{-1/(2H)}} \int_y^\infty (1 + x^{2H})^{-d/2} \, dx \, dy,$$

which converges to zero because $Hd > 3/2$.

Finally, we have

$$Z^1_\varepsilon = \int_0^\infty (1+x^{2H})^{-d/2-m} \left( \varepsilon^{1/(4H)} \int_0^{T\varepsilon^{-1/(2H)}} Z_{s,s+x}\, ds \right) dx,$$

which converges in law to a normal distribution from Proposition 5. □

PROOF OF (44). Set

$$\Phi_\varepsilon = (\log(1/\varepsilon))^{-2} \int_\mathcal{D} \prod_{i=1}^4 (\varepsilon + |t_i - s_i|^{2H})^{-d/2-m}$$
$$\times \mu^p_{1,2}(s,t) \mu^p_{3,4}(s,t) \mu^{2m-p}_{1,3}(s,t) \mu^{2m-p}_{2,4}(s,t)\, ds\, dt.$$

As in the proof of (43), it suffices to show that $\Phi_\varepsilon$ tends to zero as $\varepsilon$ tends to zero. We can also write

$$\Phi_\varepsilon = (\log(1/\varepsilon))^{-2} \mathbb{E}\left( \prod_{i,j} X^{i,j}_\varepsilon \right),$$

where the $X^{i,j}_\varepsilon$ have been defined in step 2 of the proof of (43).

From (41) it follows that the variances of the random variables $(\log(1/\varepsilon))^{-1/2} \times X^{i,j}_\varepsilon$ converge to a constant. In fact,

$$\mathbb{E}[(\log(1/\varepsilon))^{-1}(X^{i,j}_\varepsilon)^2] = \frac{(2\pi)^d 2^{2m}}{\alpha_m} (\log(1/\varepsilon))^{-1} \mathbb{E}[(I_{2m}(f_{2m,\varepsilon}))^2].$$

Then, as in the proof of part (i), it suffices to show that the family of random variables $(\log(1/\varepsilon))^{-1/2} X^{i,j}_\varepsilon$ converges in law to a normal random variable as $\varepsilon$ tends to zero. This family has the same distribution as $\{Y_\varepsilon, \varepsilon > 0\}$, where

$$Y_\varepsilon = (\log(1/\varepsilon))^{-1/2} \int_0^T \int_0^t (\varepsilon + |t-s|^{2H})^{-d/2-m} \prod_{i=1}^{2m} (B^{H,i}_t - B^{H,i}_s)\, ds\, dt,$$

and $\{B^H_t, t \geq 0\}$ is a $(2m)$-dimensional fractional Brownian motion with Hurst parameter $H$. As in the proof of part (i) we make the decomposition

$$Y_\varepsilon = Y^1_\varepsilon - Y^2_\varepsilon,$$

where

$$Y^2_\varepsilon = (\log(1/\varepsilon))^{-1/2} \varepsilon^{1/(4H)} \int_{s<T\varepsilon^{-1/(2H)}<t} (1+|t-s|^{2H})^{-d/2-m} Z_{s,t}\, ds\, dt,$$

$$Y^1_\varepsilon = (\log(1/\varepsilon))^{-1/2} \varepsilon^{1/(4H)} \int_{s<t\wedge T\varepsilon^{-1/(2H)}} (1+|t-s|^{2H})^{-d/2-m} Z_{s,t}\, ds\, dt,$$

and $Z_{s,t}$ has been defined in (16). As in the proof of part (i) we obtain

$$\|Y_\varepsilon^2\|_2 \leq (\log(1/\varepsilon))^{-1/2} \varepsilon^{1/(4H)} \int_0^{T\varepsilon^{-1/(2H)}} \int_y^\infty (1+x^{2H})^{-d/2} \, dx \, dy,$$

which clearly converges to zero. Finally, the family

$$Y_\varepsilon^1 = (\log(1/\varepsilon))^{-1/2} \int_0^\infty (1+x^{2H})^{-d/2-m} \left( \varepsilon^{1/(4H)} \int_0^{T\varepsilon^{-1/(2H)}} Z_{s,s+x} \, ds \right) dx$$

converges to a normal distribution from Proposition 6. □

REMARK 3. The limit variance appearing in (40) can be also expressed as

$$\sigma^2 = \frac{2}{(2\pi)^d} \int_{\mathbb{R}_+^3} \left[ \left( (1+x^{2H})(1+y^{2H}) \right. \right.$$

$$- \frac{1}{4}(|z+y|^{2H} + |z-x|^{2H}$$

$$\left. - |z|^{2H} - |z+y-x|^{2H})^2 \right)^{-d/2}$$

$$\left. - ((1+x^{2H})(1+y^{2H}))^{-d/2} \right] dx \, dy \, dz.$$

We do not know a simple expression for this variance, even in the case $H = 1/2$.

REMARK 4. In the particular case $H = 1/2$, and $d = 3$, the limit variance appearing in (42) is $\sigma^2 = (2\pi^2)^{-1}$. In fact, in this case we obtain

$$\sigma^2 = \frac{1}{4\pi^3} \int_{0<a+b<1} [(ab+(a+b)(1-a-b))^{-3/2}$$

$$- ((a+b)(1-a))^{-3/2} + (b(1-b))^{-3/2} - b^{-3/2}] \, da \, db.$$

Making the change of variables $x = a/(a+b)$, $y = 1/(a+b)$, it is not difficult to check that the first integral equals $2\pi$. The remaining terms are equal to $\int_0^1 b^{-3/2}(1+b-(1-b)^{-1/2}) \, db$, and with the change of variable $\sqrt{1-b} = \cos\theta$, it is clear that this integral vanishes.

**6. Technical lemmas.** We will denote by $k$ a generic constant which may be different from one formula to another one. We will decompose the region $\mathcal{T}$ defined in (9) as follows:

(46) $$\mathcal{T} \cap \{s < s'\} = \mathcal{T}_1 \cup \mathcal{T}_2 \cup \mathcal{T}_3,$$

where
$$\mathcal{T}_1 = \{(t,s,t',s'): 0 < s < s' < t < t' < T\},$$
$$\mathcal{T}_2 = \{(t,s,t',s'): 0 < s < s' < t' < t < T\},$$
$$\mathcal{T}_3 = \{(t,s,t',s'): 0 < s < t < s' < t' < T\}.$$

We will make use of the following notation:

(i) If $(t,s,t',s') \in \mathcal{T}_1$, we put $a = s' - s$, $b = t - s'$ and $c = t' - t$. On this region, the functions $\lambda$, $\rho$ and $\mu$ defined in (10) and (11) take the following values:

(47) $\quad \lambda = \lambda_1 := \lambda_1(a,b,c) := (a+b)^{2H}, \qquad \rho = \rho_1 := (b+c)^{2H},$

(48) $\quad \mu = \mu_1 := \mu_1(a,b,c) := \frac{1}{2}[(a+b+c)^{2H} + b^{2H} - c^{2H} - a^{2H}].$

(ii) If $(t,s,t',s') \in \mathcal{T}_2$, we put $a = s' - s$, $b = t' - s'$ and $c = t - t'$. On this region we will have

(49) $\quad\quad\quad \lambda = \lambda_2 := b^{2H}, \qquad \rho = \rho_2 := (a+b+c)^{2H},$

(50) $\quad\quad\quad \mu = \mu_2 := \frac{1}{2}[(b+c)^{2H} + (a+b)^{2H} - c^{2H} - a^{2H}].$

(iii) If $(t,s,t',s') \in \mathcal{T}_3$, we put $a = t - s$, $b = s' - t$ and $c = t' - s'$. On this region we will have

(51) $\quad \lambda = \lambda_3 := a^{2H}, \qquad \rho = \rho_3 := c^{2H},$

(52) $\quad \mu = \mu_3 := \frac{1}{2}[(a+b+c)^{2H} + b^{2H} - (b+c)^{2H} - (a+b)^{2H}].$

For $i = 1, 2, 3$ we set $\delta_i = \lambda_i \rho_i - \mu_i^2$, $\Theta_i = \delta_i^{-d/2} - (\lambda_i \rho_i)^{-d/2}$. Note that $\lambda_i$, $\rho_i$, $\mu_i$ and so on, $i = 1, 2, 3$, are functions of $a$, $b$ and $c$.

The following lower bounds for the determinant of the covariance matrix of $B_t^{H,1} - B_s^{H,1}$ and $B_{t'}^{H,1} - B_{s'}^{H,1}$, were obtained by Hu (2001) using the local nondeterminism property of the fractional Brownian motion [see Berman (1973)].

LEMMA 9. *There exists a constant $k$ such that, for all $a, b, c > 0$:*

(i)

(53) $\quad\quad\quad\quad\quad \delta_1 \geq k[(a+b)^{2H} c^{2H} + (b+c)^{2H} a^{2H}].$

(ii) *For $i = 2, 3$*

(54) $\quad\quad\quad\quad\quad\quad\quad\quad \delta_i \geq k \lambda_i \rho_i.$

The following lemma provides some useful inequalities.

LEMMA 10. *We have*

(55) $$k(a+b+c)^{2H-2}ac \le \mu_3 \le kb^{2H-2}ac.$$

*For $i=2,3$ we have*

(56) $$\Theta_i \le k\mu_i^2(\lambda_i\rho_i)^{-d/2-1}$$

*and*

(57) $$\Theta_i \le k(\lambda_i\rho_i)^{-d/2}.$$

PROOF. The inequalities in (55) follow from
$$\mu_3 = \tfrac{1}{2}[(a+b+c)^{2H} + b^{2H} - (a+b)^{2H} - (b+c)^{2H}]$$
$$= H(2H-1)ac \int_0^1 \int_0^1 (b+vc+ua)^{2H-2}\, du\, dv.$$

We have, for $i=2,3$,
$$\Theta_i = \left[\left(1 - \frac{\mu_i^2}{\lambda_i\rho_i}\right)^{-d/2} - 1\right](\lambda_i\rho_i)^{-d/2}.$$

The estimate (54), assuming $k<1$, implies $\frac{\mu_i^2}{\lambda_i\rho_i} \le 1-k$ and (56) holds. Moreover, (54) also implies (57). □

The following lemma is the basic ingredient in the proof of Theorem 1.

LEMMA 11. *Let $\Xi_T$ be defined by (15). Then $\Xi_T < \infty$ if and only if $dH < 3/2$.*

PROOF.

*Step* 1. Suppose $dH < 3/2$. We claim that

(58) $$\int_{[0,T]^3} \Theta_i\, da\, db\, dc < \infty$$

for $i=1,2,3$. From (53) we deduce

(59) $$\delta_1 \ge k(a+b)^H(b+c)^H a^H c^H$$
$$\ge k(abc)^{4H/3}.$$

Then, (59) together with the estimate
$$\lambda_1\rho_1 = (a+b)^{2H}(b+c)^{2H} \ge (abc)^{4H/3}$$
implies (58) for $i=1$.

To handle the case $i = 2$, we decompose the integral over the regions $\{b \geq \eta a\}$, $\{b \geq \eta c\}$ and $\{b < \eta a, b < \eta c\}$, for some fixed but arbitrary $\eta > 0$. We have, using (57),

$$\Xi := \int_{b \geq \eta a} \Theta_2 \, da \, db \, dc \leq k \int_{b \geq \eta a} \frac{da \, db \, dc}{(a+b+c)^{dH} b^{dH}}.$$

If $dH < 1$, then this integral is finite. If $1 < dH$, then

$$\Xi \leq k \int_0^T \int_0^T \frac{da \, dc}{(a+c)^{dH}} \int_{\eta a}^T b^{-Hd} \, db$$

$$\leq k \int_0^T \int_0^T a^{-4dH/3+1} c^{-2dH/3} \, da \, dc < \infty.$$

It is also easy to show that $\Xi < \infty$ in the case $1 = dH$. The case $b \geq \eta c$ can be treated in a similar way.

To deal with the case both $b < \eta a$ and $b < \eta c$, we make use of the estimate (56) and the following upper bound for $\mu_2$:

(60)
$$\mu_2 = \tfrac{1}{2}[(a+b)^{2H} - a^{2H} + (b+c)^{2H} - c^{2H}]$$
$$\leq k(a^{2H-1} + c^{2H-1})b$$

for $\eta$ small enough. In this way we obtain

$$\Theta_2 \leq k(a^{4H-2} + c^{4H-2})(a+b+c)^{-2H-dH} b^{2-2H-dH}$$
$$\leq k[a^{(2-d/3)H} b^{dH/3} + c^{(2-d/3)H} b^{dH/3}](a+b+c)^{-2H-dH} b^{-dH}.$$

Hence,

$$\int_{b<\eta a, b \leq \eta c} \Theta_2 \, da \, db \, dc$$

$$\leq k \int_{b<\eta a, b\leq \eta c} b^{-dH}(a+b+c)^{-2H-dH}$$
$$\times [a^{(2-d/3)H} b^{dH/3} + c^{(2-d/3)H} b^{dH/3}] \, da \, db \, dc$$

$$\leq k \int_{[0,T]^3} b^{-dH}(a+b+c)^{-2H-dH} a^{(2-d/3)H} b^{dH/3} \, da \, db \, dc$$

$$\leq k \int_{[0,T]^3} b^{-2dH/3} a^{-2dH/3} c^{-2dH/3} \, da \, db \, dc,$$

which is finite if $dH < 3/2$.

To handle the case $i = 3$, we decompose the integral over the regions $\{a \geq \eta_1 b, c \geq \eta_2 b\}$, $\{a < \eta_1 b, c < \eta_2 b\}$, $\{a \geq \eta_1 b, c < \eta_2 b\}$ and $\{a < \eta_1 b, c \geq$

$\eta_2 b\}$. By symmetry it suffices to consider the first three regions. We have, using (57),

$$\int_{a\geq \eta_1 b, c\geq \eta_2 b} \Theta_3 \, da \, db \, dc \leq k \int_0^T db \int_{\eta_1 b}^T \frac{da}{a^{dH}} \int_{\eta_2 b}^T \frac{dc}{c^{dH}}$$

$$\leq k \int_0^T \frac{db}{b^{2dH-2}} < \infty.$$

Let us now suppose both $a < \eta_1 b$ and $c < \eta_2 b$. Using (55) and (56) and that $H < \frac{3}{4}$ yields

$$\Theta_3 \leq k b^{4H-4} a^{2-2H-dH} c^{2-2H-dH} \leq k a^{-2dH/3} c^{-2dH/3} b^{-2dH/3},$$

which implies that the integral over this region is finite. Finally, let us consider the case $c < \eta_1 b$ and $a \geq \eta_2 b$. If $Hd > 1$, then (57) yields

$$\Theta_3 \leq k(ac)^{-Hd}$$

which is integrable. So, we can assume $H \leq \frac{1}{d} \leq \frac{1}{2}$. Then

(61)
$$\mu_3 = \tfrac{1}{2}[(a+b+c)^{2H} - (a+b)^{2H} - (c+b)^{2H} + b^{2H}]$$
$$\leq k b^{2H-1} c,$$

if $\eta_2$ is small enough. Hence, using (56) and (61) we get

$$\Theta_3 \leq k b^{4H-2} a^{-2H-dH} c^{2-2H-dH}.$$

Consequently, if $-dH + 2H + 1 < 0$,

$$\int_{c<\eta_1 b, a\geq \eta_2 b} \Theta_3 \, da \, db \, dc$$
$$\leq k \int_{c<\eta_1 b, a\geq \eta_2 b} b^{4H-2} a^{-2H-dH} c^{2-2H-dH} \, dc \, db \, da$$
$$\leq k \int_{a\geq \eta_2 b} a^{-dH-2H} b^{-dH+2H+1} \, db \, da$$
$$\leq k \int_0^T a^{-2dH+2} \, da,$$

which is finite if $dH < 3/2$. The case $-dH + 2H + 1 \geq 0$ is easier.

*Step* 2. Assume that $H = 3/(2d)$, and let us show that $\Xi_T = \infty$. It suffices to prove that

(62) $$\frac{d}{2(2\pi)^d} \int_{\mathcal{T}} \mu^2 (\lambda \rho)^{-d/2-1} \, ds \, dt \, ds' dt' = \infty,$$

because from the identity (37) this is the second moment of the second chaos of the renormalized self-intersection local time. In order to check (62) we will show that

$$A := \int_{0<a+b+c<T} (T-a-b-c)\mu_3^2(\lambda_3\rho_3)^{-d/2-1}\, da\, db\, dc = \infty.$$

With the above notation, we have, using (55), for $\varepsilon > 0$ small enough,

$$A \geq k \int_{[0,\varepsilon]^3} (b+c+a)^{4H-4}(ac)^{2-Hd-2H}\, da\, db\, dc := B.$$

If $d=2$, we get

$$B = k \int_{[0,\varepsilon]^3} \frac{1}{(a+b+c)ac}\, da\, db\, dc = \infty.$$

For $d>2$, we have $2-Hd-2H = 1/2 - 3/d > -1$. Hence,

$$B = \frac{k}{3-4H} \int_{[0,\varepsilon]^2} [(c+a)^{4H-3} - (\varepsilon+c+a)^{4H-3}](ac)^{2-Hd-2H}\, da\, dc,$$

and

$$\int_{0<a<c<\varepsilon} (c+a)^{4H-3}(ac)^{2-Hd-2H}\, da\, dc$$

$$\geq 2^{4H-3} \int_{0<a<c<\varepsilon} a^{2-Hd-2H} c^{2H-1-Hd}\, da\, dc$$

$$\geq k \int_0^\varepsilon a^{2-2Hd}\, da = \infty$$

because $2H-1-Hd < -1$. □

LEMMA 12. *If $H = 3/(2d) < 3/4$, then for $i=1,2,3$,*

$$\int_0^1 \int_0^b \Theta_i(a, b-a, 1-b)\, da\, db < \infty.$$

PROOF. Suppose first $i=1$. Then, using (53), we obtain

$$\int_0^1 \int_0^b \delta_1^{-d/2}(a, b-a, 1-b)\, da\, db$$

$$\leq k \int_0^1 \int_0^b [b(1-b) + (1-a)a]^{-3/2}\, da\, db < \infty$$

and

$$\int_0^1 \int_0^b (\lambda_1\rho_1)^{-d/2}(a, b-a, 1-b)\, da\, db$$

$$\leq k \int_0^1 b^{-3/2}[1 - (1-b)^{-1/2}]\, db < \infty.$$

For $i = 2$ we can write, using (56),

$$\int_0^1 \int_0^b \Theta_2(a, b-a, 1-b) \, da \, db$$

$$\leq k \int_0^1 \int_0^b (b-a)^{-3/2-2H} [(1-a)^{2H} + b^{2H} - (1-b)^{2H} - a^{2H}]^2 \, da \, db$$

$$\leq k \int_0^1 \int_0^b (b-a)^{-3/2-2H} (b^{2H} - a^{2H})^2 \, da \, db.$$

Then, if $\eta$ is small enough,

$$\int_{a<b, b-a<\eta a} (b-a)^{-3/2-2H} (b^{2H} - a^{2H})^2 \, da \, db$$

$$\leq k \int_{a<b, b-a<\eta a} (b-a)^{1/2-2H} a^{4H-2} \, da \, db < \infty$$

because $H < 3/4$, and

$$\int_{a<b, b-a\geq \eta a} (b-a)^{-3/2-2H} (b^{2H} - a^{2H})^2 \, da \, db$$

$$\leq k \int_0^1 \int_0^{x/\eta} x^{-3/2-2H+2\wedge 4H} \, da \, dx$$

$$= k \int_0^1 x^{-1/2-2H+2\wedge 4H} \, dx < \infty.$$

Finally, for $i = 3$, we obtain, using (56),

$$\int_0^1 \int_0^b \Theta_3(a, b-a, 1-b) \, da \, db$$

$$\leq k \int_0^1 \int_0^b a^{-3/2-2H} (1-b)^{-3/2-2H}$$

$$\times [1 + (b-a)^{2H} - (1-a)^{2H} - b^{2H}]^2 \, da \, db < \infty.$$

In fact, if we fix $\eta_1$ and $\eta_2$ small enough, it suffices to check that the above integral is finite over the regions

$$G_1 = \{a < b, a \geq \eta_1(b-a), 1-b \geq \eta_2(b-a)\},$$
$$G_2 = \{a < b, a < \eta_1(b-a), 1-b < \eta_2(b-a)\},$$
$$G_3 = \{a < b, a \geq \eta_1(b-a), 1-b < \eta_2(b-a)\},$$
$$G_4 = \{a < b, a < \eta_1(b-a), 1-b \geq \eta_2(b-a)\}.$$

On $G_1$ we have, using (57),

$$\int_{G_1} \Theta_3(a, b-a, 1-b) \, da \, db \leq k \int_{G_1} a^{-1/2-2H} (1-b)^{-1/2-2H} \, da \, db.$$

In order to estimate this integral we make the decomposition
$$G_1 \subset (G_1 \cap \{a > \varepsilon\}) \cup (G_1 \cap \{1 - b > \varepsilon\}) \cup (G_1 \cap \{a < \varepsilon, 1 - b < \varepsilon\}).$$
On the set $G_1 \cap \{a > \varepsilon\}$ the integrand is bounded by $k(1-b)^{-1/2-2H}$, and making the change of variable $y = 1 - b$, $x = b - a$, we obtain a finite integral because $H < \frac{3}{4}$. The same approach is used to handle the integral over $G_1 \cap \{1 - b > \varepsilon\}$. Finally,
$$G_1 \cap \{a < \varepsilon, 1 - b < \varepsilon\} \subset G_1 \cap \{b - a > 1 - 2\varepsilon\},$$
and, on this set the integrand is bounded.

From (55), we obtain on $G_2$
$$1 + (b-a)^{2H} - (1-a)^{2H} - b^{2H} \leq k(b-a)^{2H-2}a(1-b).$$
Hence,
$$\int_{G_2} \Theta_3(a, b-a, 1-b)\, da\, db \leq k \int_{G_2} a^{1/2-2H}(1-b)^{1/2-2H}(b-a)^{4H-4}\, da\, db.$$
On the set $G_2 \cap \{b - a > \varepsilon\}$ this integral is clearly finite, and on the set $G_2 \cap \{b - a \leq \varepsilon\}$ we have $b - a \geq 1 - 2\varepsilon$, and again the integral is finite provided $\varepsilon < 1/2$.

Fix $0 < \varepsilon < 1 - \eta_2$. Using (56), on the set $G_3 \cap \{b - a > \varepsilon\}$ we have
$$\Theta_3(a, b-a, 1-b) \leq k a^{1/2-2H}(1-b)^{1/2-2H}\varepsilon^{4H-4}$$
which is integrable. On the set $G_3 \cap \{b - a \leq \varepsilon\}$ we have
$$a = 1 - (b-a) - (1-b) \geq 1 - \varepsilon - \eta_2 > 0.$$
Then, using the estimate
$$1 + (b-a)^{2H} - (1-a)^{2H} - b^{2H} \leq |1 - b^{2H}| + |(b-a)^{2H} - (1-a)^{2H}|$$
$$\leq k[(b-a)^{2H-1} + a^{2H-1}](1-b),$$
we obtain, on the set $G_3 \cap \{b - a \leq \varepsilon\}$,
$$\Theta_3(a, b-a, 1-b) \leq k a^{-3/2-2H}(1-b)^{1/2-2H}[(b-a)^{4H-2} + a^{4H-2}]$$
$$\leq k(1-b)^{1/2-2H}[(b-a)^{4H-2} + 1],$$
which is again integrable. $\square$

We will make use of the following notation, for $i = 1, 2, 3$:

(63) $$\hat{\lambda}_i = \lambda_i + 1, \quad \hat{\rho}_i = \rho_i + 1,$$

(64) $$\hat{\delta}_i = \hat{\lambda}_i \hat{\rho}_i - \mu_i^2,$$

(65) $$\widehat{\Theta}_i = (\hat{\delta}_i)^{-d/2} - (\hat{\lambda}_i \hat{\rho}_i)^{-d/2}.$$

LEMMA 13. *Let $d > 2$ and $3/(2d) < H < 3/4$. Then for $i = 1, 2, 3$*

(66) $$\int_{\mathbb{R}_+^3} \widehat{\Theta}_i \, da \, db \, dc < \infty.$$

PROOF. The proof will be done in several steps.

*Step* 1. Suppose first $i = 1$. Then,
$$\widehat{\Theta}_1 \leq (\hat{\delta}_1)^{-d/2} = (\delta_1 + (a+b)^{2H} + (b+c)^{2H} + 1)^{-d/2}.$$

First, $\Theta_1$ is integrable in the region $\{a, b, c \geq 1\}$ because the estimate (59) yields
$$\widehat{\Theta}_1 \leq \delta_1^{-d/2} \leq k(abc)^{-2dH/3}$$
and $dH > 3/2$. On the other hand, using (53) we obtain
$$\widehat{\Theta}_1 \leq [k(a+b)^{2H}c^{2H} + k(b+c)^{2H}a^{2H} + (a+b)^{2H} + (b+c)^{2H} + 1]^{-d/2}$$
$$\leq k[\max(a^{2H} + b^{2H}, a^{2H} + c^{2H}, b^{2H} + c^{2H}) + 1]^{-d/2},$$
which implies the integrability of $\Theta_1$ on the complement of the region $\{a, b, c \geq 1\}$.

*Step* 2. Suppose that $i = 2$. Then, using (54) we obtain
$$\widehat{\Theta}_2 \leq (\delta_2 + (a+b+c)^{2H} + b^{2H} + 1)^{-d/2}$$
$$\leq [k(a+b+c)^{2H}b^{2H} + (a+b+c)^{2H} + b^{2H} + 1]^{-d/2}.$$

On $\{a, b, c \leq 1\}$ the function $\widehat{\Theta}_2$ is clearly integrable. If one of the coordinates is bigger than 1 and the other two are smaller we use the estimate
$$\widehat{\Theta}_2 \leq k[\max(a,b,c)^{2H} + 1]^{-d/2}.$$

If $a \leq 1$, $b > 1$ and $c > 1$, or if $a > 1$, $b > 1$ and $c \leq 1$ we use the estimate
$$\widehat{\Theta}_2 \leq k[\max(ab, bc)^{2H} + 1]^{-d/2}.$$

If $a > 1$, $b > 1$ and $c > 1$ we use the estimate
$$\widehat{\Theta}_2 \leq k[(a+b+c)^{2H}b^{2H}]^{-d/2} \leq k(abc)^{-2dH/3}.$$

So it only remains to check the integrability of $\widehat{\Theta}_2$ in the region
$$G = \{(a,b,c) : b \leq 1, a, c > 1\}.$$

Using inequality (56), on the region $G$ we have
$$\widehat{\Theta}_2 \leq k\mu_2^2(\hat{\lambda}_2\hat{\rho}_2)^{-d/2}$$
$$\leq k[(1+c)^{2H} + (a+1)^{2H} - c^{2H} - a^{2H}]^2(a+c)^{-Hd-2H}.$$

Consequently,

$$\int_G \widehat{\Theta}_2 \, da \, db \, dc$$
$$\leq k \int_1^\infty \int_1^\infty [(a+1)^{2H} - a^{2H}]^2 (a+c)^{-Hd-2H} \, da \, dc$$
$$\leq k \int_1^\infty \int_1^\infty a^{4H-2}(a+c)^{-Hd-2H} \, da \, dc$$
$$\leq k \int_1^\infty a^{2H-Hd-1} \, da < \infty$$

because $d > 2$.

*Step* 3. Suppose that $i = 3$. We claim that (66) holds for $i = 3$. On $\{a \geq \eta_1 b, c \geq \eta_2 b\}$ we have using (54)

$$\widehat{\Theta}_3 \leq (\delta_2 + a^{2H} + c^{2H} + 1)^{-d/2}$$
$$\leq k(a^{2H} + 1)^{-d/2}(c^{2H} + 1)^{-d/2},$$

and as a consequence, we obtain

$$\int_{a \geq \eta_1 b, c \geq \eta_2 b} \widehat{\Theta}_3 \, da \, db \, dc$$
$$\leq k \int_1^\infty db \int_{\eta_2 b}^\infty \int_{\eta_1 b}^\infty a^{-Hd} c^{-Hd} \, da \, dc \leq k \int_1^\infty \frac{db}{b^{2dH-2}} < \infty,$$

if $dH > 3/2$.

On $\{a < \eta_1 b, c < \eta_2 b\}$, using the estimates (56) and (55), we obtain

$$\widehat{\Theta}_3 \leq k b^{4H-4}(ac)^2 (a^{2H} + 1)^{-d/2-1}(c^{2H} + 1)^{-d/2-1}.$$

If $dH > 3/2$ and $H < 3/4$, then

$$\int_{a < \eta_1 b, c < \eta_2 b} \widehat{\Theta}_3 \, da \, db \, dc$$
$$\leq k \int_{[0,\infty)^2} (a^{2H} + 1)^{-d/2-1}(c^{2H} + 1)^{-d/2-1}(ac)^2 \, da \, dc \int_{a/\eta_1 \vee c/\eta_2}^\infty b^{4H-4} \, db$$
$$\leq k \int_{a<c} (a^{2H} + 1)^{-d/2-1}(c^{2H} + 1)^{-d/2-1}(ac)^2 c^{4H-3} \, da \, dc < \infty.$$

Finally, let us consider the case $c < \eta_2 b$ and $a \geq \eta_1 b$. We have, using the estimates (61) and (56),

$$\widehat{\Theta}_3 \leq k(b^{4H-2} + (a+b)^{4H-2})c^2(a^{2H} + 1)^{-d/2-1}(c^{2H} + 1)^{-d/2-1},$$

if $\eta_1$ is small enough. If $2H < 1$ we use $(a+b)^{4H-2} \le b^{4H-2}$, and if $2H > 1$ we use

$$(a+b)^{4H-2} \le (1+\eta_2)^{4H-2} b^{4H-2}.$$

In this way we obtain

$$\widehat{\Theta}_3 \le k(b^{4H-2} + a^{4H-2}) c^2 (a^{2H}+1)^{-d/2-1} (c^{2H}+1)^{-d/2-1}.$$

By decomposing the integral over the regions $\{1 < c < \eta_2 b < \frac{\eta_2}{\eta_1} a\}$, $\{c < 1 < \eta_2 b < \frac{\eta_2}{\eta_1} a\}$, $\{c < \eta_2 b < 1 < \frac{\eta_2}{\eta_1} a\}$ and $\{c < \eta_2 b < \frac{\eta_2}{\eta_1} a < 1\}$ we easily show that the integral is finite. $\square$

Department of Mathematics  
University of Kansas  
405 Snow Hall  
Lawrence, Kansas 66045  
USA  
e-mail: hu@math.ku.edu

Facultat de Matemàtiques  
Universitat de Barcelona  
Gran Via 585  
08007 Barcelona  
Spain  
e-mail: dnualart@ub.edu